\documentclass[11pt,letter]{amsart}
\usepackage[utf8]{inputenc}
\usepackage{amsmath,amsthm,amssymb,latexsym}
\usepackage[T1]{fontenc}
\usepackage{libertine} 
\usepackage{euler} 
\usepackage{eucal} 
\usepackage{color}
\usepackage{graphicx}

\def\Zorn{Zorn}

\def\sub#1{_{_{#1}}}
\def\sup#1{^{^{#1}}}

\def\mf#1{\mathfrak{#1}}
\def\mb#1{\mathbb{#1}}
\def\mc#1{\mathcal{#1}}

\def\fmod#1{f\sup{\mc{#1}}}
\def\Rmod#1{R\sup{\mc{#1}}}
\def\cmod#1{c\sup{\mc{#1}}}
\def\mod#1#2{{#1}\sup{\mc{#2}}}

\def\forza#1{\Vdash\sub{#1}}

\def\coLim#1{\underset{#1}{\text{\rm coLim}}\ }

\def\Im{\text{\rm Im}}

\def\Cond{\Rightarrow}

\def\Iff{\Leftrightarrow}

\def\ARROW#1{
\text{\begin{picture}(35,18)(15,15)         
            \put(30,25){$\sub{#1}$}
            \put(18,18){\vector(1,0){30}}
\end{picture}}}

\newcounter{numero}
\newcommand{\Numero}{\setcounter{numero}{1}(\arabic{numero}) }
\newcommand{\numero}{\addtocounter{numero}{1}(\arabic{numero}) }

\newcounter{letra}
\newcommand{\Letra}{\medskip \setcounter{letra}{1}(\alph{letra}) }
\newcommand{\letra}{\medskip \addtocounter{letra}{1}(\alph{letra}) }

\newcounter{romnumero}
\newcommand{\Romnumero}{\setcounter{romnumero}{1}(\roman{romnumero}) }
\newcommand{\romnumero}{\addtocounter{romnumero}{1}(\roman{romnumero}) }

\newcounter{bibnumero}

\newtheorem{teo}{Theorem}[subsection]                  
\newtheorem{lema}[teo]{Lemma}
\newtheorem{prop}[teo]{Proposition}
\newtheorem{cor}[teo]{Corollary}

\def\bteo{\begin{teo}}
\def\eteo{\end{teo}}
\def\bprop{\begin{prop}}
\def\eprop{\end{prop}}
\def\bcor{\begin{cor}}
\def\ecor{\end{cor}}
\def\blema{\begin{lema}}
\def\elema{\end{lema}}

\theoremstyle{remark}                               
\newtheorem{obs}[teo]{Remark}

\def\bobs{\begin{obs} }
\def\eobs{\end{obs}}
\newenvironment{dem}{ [{\it Proof\/}]\rm\hskip3mm }{\hfill$\square$\vskip5mm}       
\def\bdem{\begin{dem}}
\def\edem{\end{dem}}

\begin{document}


\pagestyle{myheadings} \markboth{G. Padilla}{Apuntes para $q$AnIntCohom}

\author{Gabriel Padilla}
\email{\rm gipadillal@unal.edu.co}

\author{Andrés Villaveces}
\email{\rm avillavecesn@unal.edu.co}

\address{Departamento de Matem\'aticas\\ 
Universidad Nacional de Colombia.\\
AK30 Cl45 Edif.404 Bogot\'a 101001000.\\ 
Tel. +00571316500 Ext.13166}

\title{Sheaves of $G$-structures and generic $G$-models}

\begin{abstract}
		In this article we give an equivariant version for the construction of
		generic models on presheaves of structures. We deal with first order structures endowed
		with a suitable action of some fixed group, say $G$; we call them $G$-structures. 
		We show that every exact presheaf of $G$-structures $\mc{M}$ has a generic $G$-model
		$\mc{M}\sup{gen}$.
\end{abstract}

\maketitle

\section*{Foreword}\label{section foreword}
This article is inspired in both the work of Caicedo for topological
sheaves of structures \cite{caicedo} and 
the geometric study of transformation groups. It is an
attempt to further the connections between Model Theory and Geometry
already opened in the work of many authors (see for instance
Macintyre\cite{macintyre} for a survey of the connections and many
open lines of work). We establish the first results towards a Model
Theoretic analysis of geometric structures beyond sheaves. We combine
the approach of Caicedo for his Model Theory on Sheaves with the
Geometric study of group actions.
Given a group $G$, a $G$-structure in a first order language
$\mc{L}=(\mc{R},\mc{F},\mc{C})$ is a structure
$\mc{A}=\left(A,\Rmod{A},\fmod{A},\mc{C}\sup{\mc{A}}\right)$ in the
usual sense, such that $G$ acts on the universe set $A$ and the 
elements of $G$ commute with the language symbols. 
A morphism of $G$-structures is a morphism of structures which is $G$-equivariant. 
The category of $G$-structures is closed under colimits~\cite{hodges}.\\[2mm]

Given a topologic space $X$; a presheaf $\mc{M}$ of $G$-structures on
$X$ is one that maps each open nbhd
$U\subset X$ to some $G$-structure $\mc{M}\sub{U}$, and each inclusion
$U\subset V$ of open nbhds to some
morphism of $G$-structures
$\mc{M}\sub{V}\ARROW{\hskip-2mm\rho\sub{UV}}\mc{M}\sub{U}$.
In this article we generalize the Generic Model Theorem 5.2 of
\cite{caicedo} to the equivariant context.
We achieve this by using only, and as far as we can, classical
presheaf techniques, combined with semantics on sheaves and
variants thereof. We study a version of the preservation of logical
``truth'',
under colimits and germs. Given an exact presheaf of $G$-structures
$\mc{M}$ on $X$; there is always a generic $G$-model $\mc{M}\sup{gen}$.
These generalizations are obtained through a
simplification in the presentation of the pointwise forcing
relation in terms of presheaves.
There are several examples of presheaves of $G$-structures coming from
a standard $G$-space and the usual functors of algebraic topology
\cite{bredon2}; such as 
singular chains, De Rham differential forms, 
usual (co)homology, intersection (co)homology and cohomology
\cite{gmp,illinois}; etc.  
New examples of structures related to semilocal geometric information
 have recently appeared in number theory and applied physics \cite{abramsky,gendron}. \\[2mm]

Our paper has the following structure: Sections 1 and 2 provides the basic framework of sheaves,
the way we will use them here: Sheaves and $G$-structures. We provide
the basic definitions and study the
preservation of validity under colimits. 
In section 3 we deal with local semantics. The final section completes the main result of 
this short paper: the definition of local forcing, construction of equivariant generic
models and a generalization of the Generic Model Theorem to this
setting.

\section{Sheaves}\label{section haces}
Recall the notion of a sheaf on a topological space $(X,\mc{T})$
\cite{godement}. 

\subsection{}\label{subsection prehaces}\label{sssection haces categoricos}
A {\bf presheaf} on $X$ with values in a category $\mf{C}$ is a map $\mc{F}$ which assigns an object 
$\mc{F}(U)\in\mf{C}$ for each open subset $U\subset X$,  and a {\bf restriction morphism} 
$\mc{F}(V)\ARROW{\hskip-2mm \rho\sub{UV}}\mc{F}(U)$
for each inclusion of open subsets $U\subset V$; in such a way that
\begin{enumerate}
	\item $\rho\sub{UU}=1\sub{\mc{F}(U)}$ for each open subset $U$ of $X$.
	\item $\rho\sub{UV}\rho\sub{VW}=\rho\sub{UW}$ for each open subsets $U\subset V\subset W$.
\end{enumerate}
In other words, $\mc{F}$ is a contravariant functor from the topology of $X$ to $\mf{C}$.
A {\bf morphism of (pre)sheaves} $\mc{F}\ARROW{T}\mc{G}$ is a natural transformation from the functor $\mc{F}$
to the functor $\mc{G}$. If the target categories are closed under inverse limits then $T$ commutes with
the limits whenever it makes sense. \medskip

A presheaf $\mc{F}$ is a {\bf  sheaf} iff for any family of open subsets and their union
\[
		\left\{U\sub{i}\ :\ i\in\mf{J}\right\}\subset\mc{T}
		\hskip2cm
		U=\bigcup_i U\sub{i}
\]
the following conditions are satisfied:  
\begin{itemize}
		\item \underline{Coherence:} If $s,s'\in\mc{F}(U)$ are such that $\rho\sub{i}(s)=\rho\sub{i}(s')$
		coincide at $\mc{F}\left(U\sub{i}\right)$ for each $i$; then $s=s'$.
		\item \underline{Exactness:} 
		If some $s\sub{i}\in \mc{F}\left(U\sub{i}\right)$ is given for each $i\in \mf{J}$ in such a way that
		$\rho\sub{ij}\left(s\sub{i}\right) = 
		\rho\sub{ij}\left(s\sub{j}\right)$ in $\mc{F}\left(U\sub{i}\cap U\sub{j}\right)$
		for each $i,j\in\mf{J}$ such that $U\sub{i}\cap U\sub{j}\neq\emptyset$; then 
		there is some $s\in\mc{F}(U)$ such that $\rho\sub{i}(s)=s\sub{i}$ in $\mc{F}\left(U\sub{i}\right)$ 
		for all $i\in\mf{J}$.
\end{itemize}

\section{$G$-structures}\label{sssection traduccion por submersion}
For a given group $G$ we introduce the notion of a $G$-structure.

\subsection{}\label{sssection estructuras}
A  {\bf language} or {\bf signature} is a triple $\mc{L}=(\mc{F},\mc{R},\mc{C})$  of sets; a set $\mc{F}$ of {\bf function symbols}, a set $\mc{R}$ 
of {\bf relation symbols} and a set $\mc{C}$ of {\bf constant symbols}. Each $f\in \mc{F}$ (resp. $R\in\mc{R}$) has an associated positive integer $n\sub{f}$ (resp. $n\sub{R}$) called its {\bf arity}. Given a signature $\mc{L}$; a {\bf structure} $\mc{M}=\left(M,\mc{F}\sup{\mc{M}},\mc{R}\sup{\mc{M}},\mc{C}\sup{\mc{M}}\right)$ is a family
of four sets satisfying
\begin{enumerate}
	\item $M\neq\emptyset$; it is the {\bf universe} of the structure $\mc{M}$.
	\item For each $f\in \mc{F}$ there is a unique function  
	$M\sup{n\sub{f}}\ARROW{\hskip-2mm f\sup{\mc{M}}}M$ in $\mod{\mc{F}}{M}$.
	\item For each $R\in\mc{R}$ there is a unique subset  $R\sup{\mc{M}}\subset M\sup{n\sub{R}}$ in $\mod{\mc{R}}{M}$.
	\item For each $c\in\mc{C}$ there is a unique $c\sup{\mc{M}}\in M$ in $\mc{C}\sup{\mc{M}}$. 
\end{enumerate}
A {\bf morphism} of structures $\mc{M}\ARROW{\alpha}\mc{N}$ is a function $M\ARROW{\alpha}N$ 
between the respective universe sets; such that:
\begin{enumerate}
		\item $\alpha\left(f\sup{\mc{M}}\left(a\right)\right)=
		f\sup{\mc{N}}\left(\alpha\left(a\right)\right)$ for each $f\in\mc{F}$
		and $a\in M\sup{n_f}$.
		\item $\alpha\left(R\sup{\mc{M}}\right)\subset R\sup{\mc{N}}$ for each $R\in\mc{R}$.
		\item $\alpha\left(c\sup{\mc{M}}\right)=c\sup{\mc{N}}$ for each $c\in\mc{C}$.
\end{enumerate} 
An {\bf isomorphism} is a bijective morphism whose inverse map is also a morphism. 
We denote the category of structures with signature $\mc{L}$ and its morphisms with the letter $\mf{M}$.\\[2mm]
A morphism $\mc{M}\ARROW{\alpha}\mc{N}$ is said to be {\bf saturated} iff  $\alpha\sup{-1}\left(R\sup{\mc{N}}\right)\subset R\sup{\mc{M}}$ for each relation symbol $R\in\mc{R}$. An {\bf embedding} (resp. a {\bf submersion}) is an injective (resp. surjective) saturated morphism. 
Given two structures $\mc{M,N}$ such that $M\subset N$; we say that $\mc{M}$ is a {\bf substructure} of $\mc{N}$ iff the inclusion map is an embedding, in 
that case we write $\mc{M}\leq\mc{N}$. \\[2mm]
Given a saturated morphism $\mc{M}\ARROW{\alpha}\mc{N}$; there is a unique substructure
$\mc{I}(\alpha)$ of $\mc{N}$ whose universe $\Im(\alpha)$ is the image set of  $\alpha$. 
From the obvious equivalence relation on $M$ induced by $\alpha$ we also obtain a quotient model $\mc{M}/\sim$ whose universe is the quotient set $M/\sim$; 
the quotien projection $\mc{M}\ARROW{q}\mc{M}/\sim$ is a submersion and the induced arrow $\mc{M}/\sim\ARROW{\overline{\alpha}}\mc{I}(\alpha)$ 
is an isomorphism.

\subsection{}\label{ssection Lformulas}
Recall the notions of formulas and validity \cite{marker}.
A term in the language is a function symbol that can be obtained, starting from a finite set of free variables and language symbols, in a finite number of steps. 
An {\bf atomic formula} is one of the form $t(v)=s(v)$ or $t(v)\in R$ where $t(v),s(v)$ are terms and $R$ is a
relation symbol. A {\bf formula} is a finite concatenation of atomic formulas and the usual logic symbols 
$\wedge,\vee,\neg,\exists,\forall$.\\[2mm]
Given a formula $\varphi(v)$ and $a\in M\sup{n}$ we say that $\mc{M}$ {\bf models} $\varphi(v)$ in $a$, and will write $\mc{M}\models\varphi(a)$; whenever $\varphi(a)$ is true.
For each morphism $\mc{M}\ARROW{\alpha}\mc{N}$ and each formula $\varphi(v)$; we will say that $\alpha$ 
{\bf preserves the validity of $\varphi$} iff whenever $a\in M\sup{n}$ and $\varphi(a)$ makes sense, 
the following conditional holds
\[
			\mc{M}\models\varphi(a)\ \Cond\ \mc{N}\models\varphi(\alpha(a))
\]

\blema\label{lema morfismos conmutan con terminos}
			{\bf [Preservation of validity under morphisms]}\\
			\Numero Morphisms commute with terms, i.e.
			given a morphism  $\mc{M}\ARROW{\alpha}\mc{N}$ and a term $t(v)$ in $n$ free variables, we have 			
			$\alpha\left(\mod{t}{M}(a)\right)=t\sup{\mc{N}}(\alpha(a))$
			for each $a\in M\sup{n}$.\\ 
			\numero Morphisms preserve the validity of formulas without $\neg,\forall$.\\
			\numero Submersions preserve the validity of formulas without $\neg$.  \\ 
			\numero If \\
			\hbox{\hskip1cm}		\Letra $\alpha$ is an isomorphism, or\\
			\hbox{\hskip1cm}  \letra $\alpha$ is an embedding and  $\varphi\left(v\right)$ is a formula
							without $\forall,\exists$; \\			
			then $M\models\varphi(a)\ \Leftrightarrow\ N\models\varphi(\alpha(a))$.
\elema
\bdem
		By induction on formulas; see  \cite[p.11-14]{marker}. 
\edem
An {\bf elementary embedding}
is an embedding $\mc{M}\ARROW{\alpha}\mc{N}$ such that, for each $a\in M\sup{n}$ and each formula $\varphi(v)$
the following equivalence holds
\[
			\mc{M}\models\varphi(a)\ \Iff\ \mc{N}\models\varphi(\alpha(a))
\]
For instance, by Lemma \ref{lema morfismos conmutan con terminos} all isomorphisms are elementary embeddings.

\subsection{}\label{sssection G-estructuras}
Fix some group $G$. A {\bf $G$-structure} is a structure $\mc{M}$ such that
\begin{enumerate}
	\item The universe $M$ is endowed with an action $G\times M\ARROW{\Phi}M$.
	\item The action commutes with the language symbols. More precisely
	\begin{enumerate}
	 \item \underline{The set of constants is invariant}: $G\mod{\mc{C}}{A}=\mod{\mc{C}}{A}$.
	 \item \underline{Relations are invariant subsets}: $G\sup{n_R}\Rmod{A}=\Rmod{A}$ 
	 for each $R\in\mc{R}$. In other words; given $\left(x\sub1,\dots,x\sub{n_R}\right)\in\Rmod{A}$
	 and $g\sub1,\dots,g\sub{n_R}\in G$, we also have
	 $\left(g\sub1 x\sub1,\dots,g\sub{n_R} x\sub{n_R}\right)\in\Rmod{A}$. 
	 \item \underline{ functions are $G$-equivariant}: For each $f\in\mc{F}$ with arity $n$,
	$g\in G$ and $x\sub1,\dots,x\sub{n}\in A$;
	\[
      \fmod{A}(g x\sub1,\dots,g x\sub{n})=
      g\fmod{A}(x\sub1,\dots,x\sub{n})
	\]
	 
	 	\end{enumerate}
\end{enumerate}

\subsection{}\label{ssection examples} Here there are some examples: \\[2mm]
\Numero  Let $G$ be a compact group 
and $X$ a topologic $G$-space. Then each $g\in G$ provides, by left multiplication, a homeomorphism $X\ARROW{g}X$. The composition with these arrows 
induces a family of chain isomorphism on singular chains
\[SC\sub{*}(X)\ARROW{g}SC\sub{*}(X).\]
 In is easy to check that $G$ acts linearly on the homology groups
$H\sub{*}(X)$, i.e. so $g(\xi+\theta)=g(\xi)+g(\theta)$ for any homology classes $\xi,\theta$; so
$H\sub{*}(X)$ is a $G$-structure. When $G$ is a Lie group and $X$ is a smooth manifold; 
a similar construction can be done for the De Rham cohomology. This can also be extended to more complicated 
(co)homology theories.\\[2mm]
\numero A countable polyhedron can be seen as a structure $\mc{M}=(\mc{S},\subset)$ in the language of posets 
$\mc{L}=(<)$, where $\mc{S}$ is a family of finite subsets in $\aleph\sub{0}$ which is hereditary for subsets (see \cite{mijares}):
if $u\subset v\in\mc{S}$ then $u\in\mc{S}$. A geometric realization $K$ of $\mc{M}$ can be obtained by setting
each $v\in\mc{S}$ to be the set of verticecs of a face in $K$. 
An easy exercise is trying to find, for a given $\mc{M}$, the biggest subgroup $G$ of the countable-order symmetric group $S\sub{\omega}$ such that $\mc{M}$ is a $G$-structure. 
For instance; the standard $n$-simplex corresponds
to the structure $\Delta\sup{n}=\left(n,\mc{P}(n)\right)$ where we take the usual identification $0=\emptyset$, and 
$n=\{0,\dots,n-1\}$ for $n>0$ and $\mc{P}(n)$ is the set parts of $n$. For this finite polyhedron, the unique
subgroup $G\leq S\sub{n}$ such that $\Delta\sup{n}$ is a $G$-structure is the trivial group $G=\{e\}$. 
On the other hand; for
$\mc{M}=\mc{P}(n)\backslash\{n\}$, which is the closure of 
$\text{max}\left(\mc{P}(n)\backslash\{n\}\right)$ by subsets (and corresponds to the boundary 
$\partial\left(\Delta\sup{n}\right)$, we have that $\mc{M}$ is a $G$-structure for any
subgroup of $S\sub{n}$. This is true because an element $g\in G$ is a permutation in the set of
vertices of $\Delta\sup{n}$; since $g$ is a bijection it preserves cardinality, so
the {\it "dimension"} on the faces is not changed by $g$.\\[2mm]
\numero The orbit set $M/G$ of a $G$-structure $\mc{M}$ is not necessarily a structure in the same 
language. For this to happen, we should consider a notion of {\it "strong"} $G$-structure
where condition 2(c) is replaced by asking each function $\fmod{M}$ to be coordinatewise equivariant;
i. e. $\fmod{M}\left(g\sub1 x\sub1,\dots,g\sub{n} x\sub{n}\right)=g\sub{1}\cdots g\sub{n}\fmod{M}(x\sub1,\dots,x\sub{n})$;
for $g\sub{1},\dots,g\sub{n}\in G$; $x\sub{1},\dots,x\sub{n}\in M$. Then the quotient function
induced by $\fmod{M}$ makes sense in $\left(\mc{M}/G\right)\sup{n}$. A similar work to the one written here
can be carried out for orbit structures coming from strong $G$-structures: The category is closed by
colimits, there is an open semantic and one can prove the existence of generic orbit models. More over, 
$\mc{M}\sup{gen}/G=\left(\mc{M}/G\right)\sup{gen}$ (see next sections for the definition of the
generic model $\mc{M}\sup{gen}$). Finite polyhedra in example (2) are strong $G$-structures (for a suitable
choice of $G$). Singular chains, smooth forms and (co)homology classes are {\it not} examples of
strong $G$-structures.

\subsection{}\label{ssection orbits and stabilizers}
A {\bf morphism} of $G$-structures is a $G$-equivariant morphism of structures. 
By a {\bf $G$-substructure} of $\mc{M}$ we mean a $G$-invariant substructure. 
The composition of $G$-equivariant morphisms 
(resp. embeddings, submersions, elementary embeddings) is again an
arrow of this kind. 
The family of $G$-structures and $G$-equivariant morphisms
(resp. embeddings, submersions, elementary embeddings) 
is a category, we will denote it by $\mf{M}\sub{G}$ (resp. $\mf{M}\sup{\leq}\sub{G}$, 
$\mf{M}\sup{\geq}\sub{G}$, $\mf{M}\sup{\prec}\sub{G}$).

\subsection{}\label{ssection notational convention}\label{ssection limites de estructuras}
In what follows we will study some properties of $G$-structures, orbit structures and their limits. 
We will assume the following conventions: For an inverse system of $G$-structures
$\{\mc{M}\sub{i}:i\in\mc{D}\}$ and $a\in M\sub{i}$ for some $i\in \mc{D}$ we write $[a]$ for the
germ of $a$ in the colimit $\mc{M}=\coLim{i\in\mc{D}}\mc{M}\sub{i}$. As we will show in this \S;
there is a well defined germ action of $G$ on $\mc{M}$. We write $\langle a\rangle$ for the orbit
of $[a]$ in the colimit structure.

\bprop\label{prop limites de estructuras}
	$\mf{M}\sub{G}$, $\mf{M}\sub{G}\sup{\leq}$, $\mf{M}\sub{G}\sup{\prec}$ and $\mf{M}\sub{G}\sup{\geq}$ are 
	closed by colimits.
\eprop
\bdem 
		We proceed by steps:\\[2mm]
		$\bullet$ \underline{Definition of a colimit}: See \cite[p.49-52]{hodges}.
		Given an inverse system of structures $(\mc{D},\leq)\ARROW{}\mc{M}$; 
		the colimit structure $\mc{M}=\coLim{i\in\mc{D}}\mc{M}\sub{i}$ always exists. 
		Its universe set $M$ is the quotient of the disjoint union
		$\underset{i}{\sqcup}M\sub{i}$ by the equivalence
		relation 
		{\small\[
			\forall i\ \forall j\ \forall x\in M\sub{i}\ \forall y\in M\sub{j} 
			\left(x\sim y \Leftrightarrow\ \exists k\leq i,j: (\rho\sub{ki}(x)=\rho\sub{kj}(y))\right) 
		\]}
		Denote $[x]$ the equivalence class of $x\in M\sub{i}$.
		For $f\in\mc{F}$ define $\mod{f}{M}([x])=\left[f\sup{\mc{M}_i}(x)\right]$.
		Also, given $R\in\mc{R}$ let $\mod{R}{M}=\coLim{i\in\mc{D}}{R\sup{\mc{M}_i}}$. Finally, if $c\in \mc{C}$ 
		take $\cmod{M}=\left[c\sup{\mc{M}_i}\right]$ for any
		$i\in\mc{D}$. Notice that the quotient maps $\mc{M}\sub{i}\ARROW{q_i}\mc{M}$ are morphisms. \\[2mm]
		$\bullet$ \underline{$\mf{M}\sub{G}$ is closed by limits}: If the above is a directed system of $G$-structures
		with equivariant arrows; then there is a $G$-action  on the universe $M$ of $\mc{M}$, given by  
		\[
					g[x]=[gx]
		\]
		and each quotient map $q\sub{i}$ is $G$-equivariant. Next we show that:\\
		\begin{itemize}
		\item[\Letra] Functions are equivariant: 
		Let $f\in\mc{F}$ be a function symbol with arity $n$. For each $x\sub{1},\dots,x\sub{n}$ in $M\sub{i}$ 
		and $g\in G$, 
		\[
			      \begin{array}{ll}
				    \fmod{M}(g[x\sub{1}],\dots,g[x\sub{n}]) & =\fmod{M}([gx\sub{1}],\dots,[gx\sub{n}])
				    =\left[\fmod{M}\sub{i}(gx\sub{1},\dots,gx\sub{n})\right]\\[2mm]
				    &=\left[g\fmod{M}\sub{i}(x\sub{1},\dots,x\sub{n})\right]=g\left[\fmod{M}\sub{i}(x\sub{1},
				    \dots,x\sub{n})\right]\\[2mm]
				    &=g\fmod{M}([x\sub{1}],\dots,[x\sub{n}])
			      \end{array}
		\]
		\item[\letra] Relations are invariant: Consider the germ $[x\sub{1},\dots,x\sub{n}]\in\Rmod{M}$ of a tuple
		$(x\sub{1},\dots,x\sub{n})\in R\sup{\mc{M}_i}$, and $g\sub{1},\dots, g\sub{n}\in G$. Then, by our assumption,
		$(g\sub{1}x\sub{1},\dots,g\sub{n}x\sub{n})\in R\sup{\mc{M}_i}$. Since the restriction maps are
		equivariant, 
		{\small\[
			(\rho\sub{ij}(g\sub{1}x\sub{1}),\dots,\rho\sub{ij}(g\sub{n}x\sub{n}))			
			=(g\sub{1}\rho\sub{ij}(x\sub{1}),\dots,g\sub{n}\rho\sub{ij}(x\sub{n}))
			\in \rho\sub{ij}\sup{n}\left(\mc{R}\sup{M_i}\right)
			\subset R\sup{\mc{M}_j}
		\]}
		for any $j>i$. This implies 
		$[g\sub{1}x\sub{1},\dots,g\sub{n}x\sub{n}]\in \Rmod{M}$.
		\item[\letra] The set of constants if invariant: For $c\in\mc{C}$, 
		$g\left[\cmod{M}\right]=\left[gc\sup{\mc{M}_i}\right]$ for some (any) $i$. Since
		$gc\sup{\mc{M}_i}\in\mc{C}\sup{\mc{M}_i}$, we are done. 
		\end{itemize}		 
		$\bullet$ \underline{$\mf{M}\sub{G}\sup{\leq}$, $\mf{M}\sub{G}\sup{\prec}$ and $\mf{M}\sub{G}\sup{\geq}$ 
		are closed by limits}:
		Since we put no condition on the restriction morphisms, the other subcategories are still closed
		under inverse limits.  
\edem

\bprop\label{cor limites de estructuras}
	Let $\mc{M}=\coLim{i\in\mc{D}}\mc{M}\sub{i}$ be a colimit of $G$-structures 
	and $a\in M\sup{n}\sub{i}$ for some $i\in\mc{D}$. If 
	$\varphi(v)$ has no $\neg,\forall$; then
	$\mc{M}\models\varphi([a])$ if and only if there is some $j\leq i$ such that 
	$\mc{M}\sub{j}\models\varphi(\rho\sub{ji}(a))$.
\eprop
\bdem Let's show the implication $(\Cond)$. If $\mc{M}\models\varphi([a])$; then we can check, by induction, 
the following cases:
	\begin{itemize}
			\item[\Letra] \underline{$\varphi$ is $t(v)=s(v)$:}
			If $\mc{M}\models\varphi([a])$ then 
			$t\sup{\mc{M}}([a])=s\sup{\mc{M}}([a])$
			so $\left[t\sup{\mc{M}\sub{i}}(a)\right]=\left[s\sup{\mc{M}\sub{i}}(a)\right]$.
			There is some $j\leq i$ such that 
			$\rho\sub{ji}\left(t\sup{\mc{M}\sub{i}}(a)\right)=
			\rho\sub{ji}\left(s\sup{\mc{M}\sub{i}}(a)\right)$. Since $\rho\sub{ji}$
			is a morphism, by \S\ref{lema morfismos conmutan con terminos}-(1) we get  
			$t\sup{\mc{M}\sub{j}}(\rho\sub{ji}(a))=s\sup{\mc{M}\sub{j}}(\rho\sub{ji}(a))$;
			so $\mc{M}\sub{j}\models\varphi(\rho\sub{ji}(a))$.\\[1mm]

			\item[\letra] \underline{$\varphi(v)$ is $t(v)\in R$:} 
			$\mc{M}\models\varphi([a])$ iff 		
			$t\sup{\mc{M}}([a])=\left[t\sup{\mc{M}\sub{i}}(a)\right]\in\Rmod{M}$.
			By step (4) in the proof of \S\ref{prop limites de estructuras}
			$\exists j\leq i$ such that 
			$t\sup{\mc{M}\sub{j}}(\rho\sub{ji}(a))=
			\rho\sub{ji}\left(t\sup{\mc{M}\sub{i}}(a)\right)\in R\sup{\mc{M}\sub{j}}$, so 
			$\mc{M}\sub{j}\models\varphi(\rho\sub{ji}(a))$. \\[1mm] 
		\end{itemize} 
		So it holds for atomic formulas. Next we apply the inductive hypothesis for
		\begin{itemize} 
			
			\item[\letra] \underline{$\varphi(v)$ is $\psi(v)\wedge\theta(v)$:} 
			By induction, assume the statement for both
                        $\psi(v)$ and $\theta(v)$. 
			$\mc{M}\models\varphi([a])$ iff $\mc{M}\models\psi([a])$ and  $\mc{M}\models\theta([a])$.
			Then $\exists k,k'\leq i$ such that 
			$\mc{M}\sub{k}\models\psi(\rho\sub{ki}(a))$ and 
			$\mc{M}\sub{k'}\models\theta(\rho\sub{k'i}(a))$.
			Take any $j\leq k,k'$ and 
			$\rho\sub{ji}(a)=\rho\sub{k'j}(\rho\sub{k'i}(a))=
			\rho\sub{kj}(\rho\sub{ki}(a))$. By \S\ref{lema morfismos conmutan con terminos}-(4.b)
			we get $\mc{M}\sub{j}\models\psi(\rho\sub{ji}(a))$
			and $\mc{M}\sub{j}\models\theta(\rho\sub{ji}(a))$. Therefore, 
			$\mc{M}\sub{j}\models\varphi(\rho\sub{ji}(a))$.\\[1mm] 
			\item[\letra] \underline{$\varphi(v)$ is $\psi(v)\vee\theta(v)$:} 
			This is similar to the previous step.
			\item[\letra] \underline{$\varphi(v)$ is $\exists w\psi(v,w)$:} 
			$\mc{M}\models\varphi([a])$ iff there is some germ
			$[b]$ such that $\mc{M}\models\psi([a],[b])$. By induction we can assume that
			there is some $j\leq i$ such that $b\in M\sub{j}\sup{n'}$ and 
			$\mc{M}\sub{j}\models\psi(\rho\sub{ij}(a),b)$. Then
			$\mc{M}\sub{j}\models\varphi(\rho\sub{ij}(a))$.
		\end{itemize} 
	This proves $\Cond$; the converse hods by \S\ref{lema morfismos conmutan con terminos}-(2). 
\edem

\section{Local semantics}\label{section local semantics}
Next we recall the notions of point and open semantics on a presheaf of structures. These are natural continuous extensions 
of Propositions \ref{cor limites de estructuras}.

\subsection{}\label{ssection haces de estructuras}\label{prehaces de estructuras}
A {\bf presheaf of $G$-structures} on $X$ is a presheaf $\mc{T}\ARROW{\mc{M}}\mf{M}\sub{G}$
in the sense of \S\ref{subsection prehaces}. Each open subset $U$ of  $X$ is sent to some $G$-structure $\mc{M}\sub{U}$ 
and each inclusion of open subsets $U\subset V$ is mapped to the corresponding equivariant restriction morphism
$\mc{M}\sub{V}\ARROW{\hskip-2mm \rho\sub{UV}}\mc{M}\sub{U}$. When the group $G$ is trivial we obtain,
in particular, a presheaf of structures as usual \cite{caicedo}.

\subsection{Point semantics}\label{ssection point forcing}
Fix some presheaf of $G$-structures $\mc{M}$ on $X$ and a point $x\in X$. 
Let $\varphi(v)$ be a formula in free variables $v=(v\sub1,\dots,v\sub{n})$. Given an open nbhd $U\ni x$
and some element $a\in\mc{M}\sub{U}$; we say that $\mc{M}$ {\bf forces}
$\varphi(a)$ at $x$, and we will write 
\[
      \mc{M}\forza{x}\varphi(a)
\]
in the following cases:\\[2mm]
\Numero \underline{$\varphi(v)$ has no $\neg,\forall$}: Iff there is some open nbhd $U\supset V\ni x$ 
such that $\mc{M}\sub{V}\models\varphi\left(\rho\sub{UV}(a)\right)$. Notice that, by Proposition \ref{cor limites de estructuras}
this is equivalent to require that $\mc{M}\sub{x}\models\varphi\left([a]\sub{x}\right)$.\\
\numero \underline{$\varphi(v)$ is $\neg\psi(v)$}: Iff there is some
open nbhd $U\supset V\ni x$ such that $\mc{M}\not\Vdash_{y}\psi(a)$ for
all $y\in V$.\\
\numero \underline{$\varphi(v)$ is $\psi(v)\rightarrow\nu(v)$}: Iff there is some open nbhd $U\supset V\ni x$ such that, for all $y\in V$; if
		$\mc{M}\forza{y}\psi(a)$ then $\mc{M}\forza{y}\nu(a)$.\\
\numero \underline{$\varphi(v)$ is $\forall w\psi(v,w)$}: Iff  there is some open nbhd $U\supset V\ni x$ such that,
for each $y\in V$ and each $b\in\mc{M}\sub{V}$, we have $\mc{M}\forza{y}\psi\left(a,b\right)$.

\bprop\label{lema forcing for presheaves}
      On categoric sheaves of structures; the above definition of point semantics is equivalent 
      to the one provided at \cite{caicedo} for topologic sheaves.
\eprop
\bdem
      Let $\mc{M}$ be a sheaf of structures, $\xi=(E,p,X)$ the topologic sheaf induced by $\mc{M}$; this is
      a sheaf of structures as defined at \cite{caicedo}. We  must show that
      \begin{equation}\label{eq equivalence point forcing}
	    \mc{M}\forza{x}\varphi(a)\ \Iff\ \xi\forza{x}\varphi(\sigma)
      \end{equation}
      for some local section $\sigma$ defined at $x$. We will do this by induction on $\varphi(v)$.
    Let $U\subset X$ be an open set. By \cite[p.110]{godement}, $\mc{M}\sub{U}$ is the structure of local sections of $\xi$ defined at $U$.
    Each local section $\sigma$ defined at $U$ is given in terms of some element $a\in M\sub{U}$ by the canonic representation map
    which sends each point $y\in U$ to the germ of $a$ in the colimit structure $\mc{M}\sub{y}$; we will write this situation with the 
    identity 
    \[
	  \sigma(y)=[a]\sub{y}\hskip1cm \forall y\in U
    \]
    If $\varphi(v)$ is atomic then, by Proposition \ref{cor limites de estructuras} and Definition 3.1-(1) 
    at \cite[p.15]{caicedo},
    \[
      \mc{M}\forza{x}\varphi(a)\ \Iff\ \mc{M}\sub{x}\models\varphi([a]\sub{x})\ 
      \Iff\ \mc{M}\sub{x}\models\varphi(\sigma(x))\ 
      \Iff\ \xi\forza{x}\varphi(\sigma)
    \]
    as desired. Suppose that $\varphi(v)$ is $\alpha(v)\wedge\beta(v)$ and the statement holds for 
    $\alpha,\beta$. If $\xi\forza{x}\varphi(\sigma)$ then, by the corresponding definition, 
    $\xi\forza{x}\alpha(\sigma)$ and $\xi\forza{x}\beta(\sigma)$. By induction, 
    $\mc{M}\forza{x}\alpha(a)$ and $\mc{M}\forza{x}\beta(a)$. Let $V\sub{\alpha},V\sub{\beta}\subset U$
    be open nbhds of $x$ such that $\mc{M}\sub{V_\alpha}\models\alpha\left(\rho\sub{UV_\alpha}(a)\right)$
    and $\mc{M}\sub{V_\beta}\models\beta\left(\rho\sub{UV_\beta}(a)\right)$. Take 
    $V=V\sub{\alpha}\cap V\sub{\beta}$. Then, by Proposition \ref{lema morfismos conmutan con terminos}-(2);
    $\mc{M}\sub{V}\models\varphi\left(\rho\sub{UV}(a)\right)$; so $\mc{M}\forza{x}\varphi(a)$. This proves
    one implication, the converse is straightforward. If $\varphi(v)$ is $\alpha(v)\vee\beta(v)$ one
    can proceed in a similar way. The case when $\varphi(v)$ is $\exists w\psi(v,w)$ is straightforward. 
    This proves the equivalence (\ref{eq equivalence point forcing}) above
    for any formula $\varphi(v)$ without $\neg,\forall$.  
    Conditions \S\ref{ssection point forcing} (2),...,(5) are equivalent to their corresponding statements
    at Definition 3.1 (4),...,(7) in \cite[p.16]{caicedo}. This finishes the proof.
\edem

\subsubsection{}\label{ssection properties point forcing}
Given a point $x\in X$, an open nbhd $U\ni x$, a formula $\varphi(v)$ in free variables 
$v=(v\sub1,\dots,v\sub{n})$ and some $a\in M\sub{U}\sup{n}$; the following properties hold:  
\begin{itemize}
	\item[\Letra] \underline{Local Semantics}: $\mc{M}\forza{x}\varphi(a)$ iff there is some open nbhd 
	$U\supset V\ni x$  such that $\mc{M}\forza{y}\varphi(a)$ for all $y\in V$.
	\item[\letra] \underline{Classical Semantics}: For  an
          isolated point $x\in X$ we get 
	$\mc{M}\forza{x}\varphi(a) \Iff \mc{M}\sub{x}\models\varphi([a]\sub{x})$.			
	\item[\letra] \underline{Excluded Middle Principle}: 
	$\mc{M}\forza{x}\forall u\forall v(u=v)\vee(u\neq v)$
	iff \Romnumero $x$ has an open nbhd $U_1$ such that, for all $a,b\in\mc{M}\sub{U_1}$,
	$a=b$. Or, \romnumero $x$ has an open nbhd $U_2$ such that, for all $a,b\in\mc{M}\sub{U_2}$,
	$a\neq b$. When $\mc{M}$ is a sheaf and $X$ is Hausdorff, this means that the induced topologic 
	sheaf is also Hausdorff in some nbhd of $x$. 
\end{itemize}

\subsection{Open semantics}\label{sssection local forcing}
Given a presheaf of $G$-structures $\mc{M}$, an open nbhd $U\subset X$ and some 
$a\in\mc{M}\sub{U}$; we say that $\mc{M}$ {\bf forces $\varphi(a)$ in $U$}, and we write
$\mc{M}\forza{U}\varphi(a)$, iff $\mc{M}\forza{x}\varphi(a)$ for all $x\in U$. 
By \S\ref{ssection properties point forcing}-(1),
\[
    \mc{M}\forza{x}\varphi(a)\ \Iff\ \text{there is some nbhd }U\supset V\ni x\text{ such that }
    \mc{M}\forza{V}\varphi(a)
\]

Also, the validity of $\varphi(a)$ is related to the topology of $X$ as follows:  
\begin{itemize}
	\item[\Letra] \underline{Restrictions}: If $U\subset V$ then $\mc{M}\forza{V}\varphi(a)\ \Cond\ \mc{M}\forza{U}\varphi(a)$. 
	\item[\letra] \underline{Coverings}: $\mc{M}\forza{U_i}\varphi\left(a|\sub{U_i}\right)\ \forall i\ \Cond 
	\mc{M}\forza{\underset{i}{\cup} U_i}\varphi(a)$.
	\item[\letra] \underline{Existencial quantifier}: $\mc{M}\forza{U}\exists\nu\varphi(a,\nu)$ iff there 
	is an open covering $\underset{i}{\cup} U_i\supset U$ and some $b\sub{i}\in\mc{M}\sub{U_i}$ for each
	$i$; such that $\mc{M}\forza{U_i} \varphi\left(a,b\sub{i}\right)$ for each $i$.		
\end{itemize}

\bprop\label{prop maximum principle}	
	{\bf [Maximum principle]} 
	Let $\mc{M}$ be an exact presheaf of $G$-structures on $X$, $x\in X$ a point, 
	$U\ni x$ an open nbhd and and $a$ in $\mc{M}\sub{U}$.
	If $\mc{M}\forza{U}\exists\nu\varphi(a,\nu)$ then there is 
	some open subset $V\subset U$ and $b\in\mc{M}\sub{V}$ such that $U\subset\overline{V}$ and
	$\mc{M}\forza{V}\varphi(a,b)$.
\eprop
\bdem
	This is a traslation of Theorem 3.3 in \cite[p.18]{caicedo}.	
	Let $\mf{X}=\underset{V\subset U}{\cup}M\sub{V}$ be the union of the universes of structures defined
	on open subsets of $U$; this is, by our assumptions, a nonempty set. 
	Consider in $\mf{X}$ the following partial order relation: For $b\in M\sub{V}$
	and $b'\in M\sub{W}$ define
	\[
	    b\leq b'\ \Iff\ V\subset W\text{ and }\rho\sub{VW}(b')=b
	\]
	If $\{a\sub{i}\in M\sub{V_i}:i\in\mf{J}\}$ is a chain in $\mf{X}$ then, since
	$\mc{M}$ is exact, for $V=\underset{i}{\cup}V\sub{i}$ there is some $a\in M\sub{V}$
	such that $a\sub{i}\leq a$ for all $i$. Since $V\subset U$ we get $a\in\mf{X}$.
	By the Zorn Lemma; there is some maximal element $a\in M\sub{W}\subset\mf{X}$.
	By the maximality of $a$, $W$ is dense in $U$.
\edem

\section{Equivariant generic models}\label{section equivariant generic models}
In the rest of this work we show how the models constructed at 
\S\ref{section equivariant generic models} are {\it generic}. Let us fix a presheaf of $G$-structures
$\mc{M}$ on $X$. 

\subsection{}\label{ssection filtros genericos}
A {\bf filter} in $X$ is a family of open subsets $\mb{F}$ which is closed by finite intersections and 
open supsets: For any $U,V$ open in $X$,
\begin{enumerate} 
		\item $U,V\in\mb{F}\ \Cond\ U\cap V\in\mb{F}$.
		\item $V\in\mb{F}$ and $U\supset V\ \Cond\ U\in\mb{F}$
\end{enumerate} 
Notice that $\mb{F}$ is trivial iff $\emptyset\in\mb{F}$.
We say that $\mb{F}$ is {\bf maximal} iff it is not properly contained in any other filter. 
A straightforward application of the \Zorn's Lemma  shows that there are maximal filters. 
A non trivial filter of open subsets 
$\mb{F}$ in $X$  is {\bf generic with respect to $\mc{M}$ } iff:
\begin{enumerate}
		\item For each formula $\varphi(v)$ in free variables $v=(v\sub1,\dots,v\sub{n})$, 
		$U\in\mb{F}$ and $a\in M\sub{U}\sup{n}$; there is some  
		$U\supset V\in\mb{F}$ such that $\mc{M}\forza{V}\varphi(\sigma)$ or $\mc{M}\forza{V}\neg\varphi(\sigma)$. 
		\item For each formula $\varphi(v,w)$ in free
                  variables $v=(v\sub1,\dots,v\sub{n})$ and
		$w=(w\sub1,\dots,w\sub{m})$; each $U\in\mb{F}$ and $a\in M\sub{U}\sup{n}$; 
		if $\mc{M}\forza{U}\exists w\varphi(a,w)$ then		
		there is some $U\supset V\in\mb{F}$ and $b\in M\sub{V}\sup{m}$ such that 
		$\mc{M}\forza{V}\varphi(a,b)$.
\end{enumerate}

\bprop\label{prop generic filters2}
      {\bf [Existence of generic filters]} If $\mc{M}$ is an exact presheaf of $G$-structures; then  
      every maximal filter in $X$ is generic with respect to $\mc{M}$.
\eprop
\bdem
      For condition \S\ref{ssection filtros genericos}-(1) apply the same proof of 
      Theorem 5.1 at \cite[p.27]{caicedo}. For condition
      \S\ref{ssection filtros genericos}-(2) the main argument is the maximum principle which,
      in our context, only requires the hypothesis of exactness.
\edem

\subsection{}\label{ssection modelos genericos}
A {\bf generic model} for a presheaf of $G$-structures $\mc{M}$ is the colimit structure
\[
    \mc{M}\sup{gen}=\coLim{U\in\mb{F}}\mc{M}\sub{U}
\]
of the structures on a generic filter $\mb{F}$ of $X$.

\subsection{}\label{ssection Godel translation}
Let us show the behavior of the forcing relation under double negations.
We start with two easy statements, the proofs are left to the reader
who can go to \cite{caicedo} for more details.

\blema\label{lemma double negation}
    Let $\varphi(v)$ be a positive formula. Then 
	$\mc{M}\forza{U}\neg(\neg\varphi(a))$ iff there is some open nbhd $V\subset U$
	such that $V$ is dense in $U$ and $\mc{M}\forza{V}\varphi(a)$.
\elema

\blema\label{lemma filters and dense subsets}
    Let $\mb{F}$ be a maximal filter of open nbhds in $X$, and $U\in \mb{F}$.
    If $V\subset U$ is open and dense in $U$, then $V\in \mb{F}$.
\elema

The {\bf G\"odel translation} $\varphi\sub{\mathbb{G}}$ of some formula $\varphi$ is defined, by induction, as follows: 
\begin{itemize}
 \item $\varphi\sub{\mathbb{G}}$ is $\neg(\neg \varphi)$ for an atomic formula $\varphi$.
 \item $(\varphi\wedge\psi)\sub{\mathbb{G}}=\varphi\sub{\mathbb{G}}\wedge\psi\sub{\mathbb{G}}$.
 \item $(\varphi\vee\psi)\sub{\mathbb{G}}=\neg
 \left(\neg\varphi\sub{\mathbb{G}}\wedge\neg\psi\sub{\mathbb{G}}\right)$.
 \item $(\neg\varphi)\sub{\mathbb{G}}=\neg\left(\varphi\sub{\mb{G}}\right)$.
  \item $(\forall v\varphi)\sub{\mb{G}}=\forall v\left(\varphi\sub{\mathbb{G}}\right)$.
 \item $(\exists v\varphi)\sub{\mb{G}}=\neg\forall v\left(\neg\varphi\sub{\mathbb{G}}\right)$. 
\end{itemize}

\bteo\label{teo genericidad}
    {\bf[Equivariant generic model theorem]} Let $\mc{M}$ be a
    sheaf of $G$-structures
    on $X$ and $\mc{M}\sup{gen}$ the generic model induced by some
    generic filter    
    $\mb{F}$ on $X$. For each formula $\varphi(v)$, 
    $U\in\mb{F}$ and $a\in\mc{M}\sub{U}$; the following statements are
    equivalent: 
    \begin{enumerate}
	  \item $\mc{M}\sup{gen}\models\varphi([a])$.
	  \item
            $\mc{M}\forza{V}\varphi\sub{\mb{G}}\left(a\right)$
            for some $U\supset V\in\mb{F}$.
	  \item $\{x\in U:
            \mc{M}\forza{x}\varphi\sub{\mb{G}}(a)\}\in\mb{F}$.
    \end{enumerate}
\eteo
\bdem We proceed by steps.\\[2mm]
      $(1)\Iff(2)$: Let us proceed by induction. 
      \begin{itemize}
	  \item \underline{$\varphi(v)$ is atomic}: 
	  By Proposition \S\ref{cor limites de estructuras},
          $\mc{M}\sup{gen}\models\varphi([a])$
	  iff there is some open nbhd $U\supset V\in\mb{F}$ such that 
	  $\mc{M}\sub{U}\models\varphi\left(\rho\sub{UV}(a)\right)$. 
	  The definition of local forcing \S\ref{sssection local forcing},
	  implies that $\mc{M}\forza{V}\varphi\left(a\right)$. 
	  By definition, $\varphi\sub{\mb{G}}=\neg(\neg\varphi)$. Finally, 
	  by Lemma \ref{lemma double negation},
	  $\mc{M}\forza{V}\varphi\sub{\mb{G}}\left(a\right)$.
          Conversely;
          if 
	  $\mc{M}\forza{V}\varphi\sub{\mb{G}}\left(a\right)$
          for some $V\in\mb{F}$
	  then there is some open nbhd $W\subset V$ such that $W$ is
          dense in $V$
	  and $\mc{M}\forza{W}\varphi\left(a\right)$. By
          Lemma \ref{lemma filters and dense subsets},
	  $W\in\mb{F}$ and, by Proposition \S\ref{cor limites de estructuras}, 
	  $\mc{M}\sup{gen}\models\varphi([a])$. 
	  \item \underline{$\varphi(v)$ is $\alpha(v)\wedge\beta(v)$}: 
	  $\mc{M}\sup{gen}\models\varphi([a])$ iff
          $\mc{M}\sup{gen}\models\alpha([a])$ and
	  $\mc{M}\sup{gen}\models\beta([a])$. By induction;
          this happens iff there are open nbhds 
	  $U\supset V\sub1,V\sub2$ in $\mb{F}$ satisfying  
	  $\mc{M}\forza{V_1}{\alpha}\sub{\mb{G}}\left(a\right)$
          and
	  $\mc{M}\forza{V_2}{\beta}\sub{\mb{G}}\left(a\right)$. By
	  inductive definition,
	  \[
		{\alpha}\sub{\mb{G}}\wedge
                {\beta}\sub{\mb{G}}=
		\left({\alpha\wedge\beta}\right)\sub{\mb{G}}=
                \varphi\sub{\mb{G}}
	  \]
	  Then $V=\left(V\sub1\cap V\sub2\right)\in\mb{F}$ and
	  $\mc{M}\forza{V}\varphi\sub{\mb{G}}\left(a\right)$,
          as desired.
	  The converse is straightforward.
	  \item \underline{$\varphi(v)$ is $\neg\alpha(v)$}: 
	  $\mc{M}\sup{gen}\models\varphi([a])$ iff 
	  $\mc{M}\sup{gen}\not\models\alpha([a])$. By induction; 
	  $\mc{M}\sub{V}\not\models\alpha\sub{\mb{G}}(a)$ for all
	  $U\supset V\in\mb{F}$, let us pick just one $V$ from these.
	  Since
          $\varphi\sub{\mb{G}}=\neg\alpha\sub{\mb{G}}$ 
	  we get that
	  $\mc{M}\sub{V}\models\varphi\sub{\mb{G}}$. 
	  \item \underline{$\varphi(v)$ is $\forall v\alpha(v)$}:
            Proceed as in the previous step.
    \end{itemize}
    The other inductive steps are easy to check; we leave them to the
    reader. Finally, the equivalence
    $(2)\Iff(3)$ can be seen as before, with an inductive proof. It
    is, essentially, a consequence
    of Lemmas \ref{lemma double negation} and \ref{lemma filters and
      dense subsets}.
\edem

\section*{Acknowledgments}
The authors would like to thank P. Zambrano, W. Páez and
X. Caicedo for some useful remarks. The second author would like to
thank Tim Gendron for his insightful comments on the
general field of application of model theoretic sheaves (and
$G$-sheaves) to problems in Geometry and Number Theory: The
study of those applications was an important motivation for this paper.


\begin{thebibliography}{99}

\bibitem{abramsky} ABRAMSKY, S.; MANSFIELD, S.  \& SOARES, R. {\sl The cohomology of Non-Locality and Contextuality}.
8th. Int. Workshop on QPL. EPTCS 95 (2012) pp. 1-14 .

\bibitem{bredon2} BREDON, G. {\sl Topology and geometry.\/} Graduate Texts in Mathematics Vol. {\bf 139}. Springer-Verlag.
Berlin (1993).

\bibitem{caicedo} CAICEDO, X. {\sl L\'ogica de loas haces de estructuras\/}. Rev. Acad. Colomb. Cienc. {\bf 19}-74. pp.569-586 (1995).

\bibitem{chang} CHANG, C. \& KEISLER, H. {\sl Model Theory.\/}
 Studies in Logic and the Foundations of Math. Vol. {\bf 73} Elseviere- Amsterdam (1993).

\bibitem{gendron} GENDRON, T. {\sl Real algebraic number theory II.\/} Arxiv.math NT.1208.4334 (2012). 
 


\bibitem{godement}  GODEMENT, R. {\sl Topologie Alg\'ebrique et Th\'eorie des
Faisceaux.\/} Hermann. Paris (1958).

\bibitem{gmp}  GORESKY, M. \& MACPHERSON, R. {\sl Intersection homology theory.\/} Topology.
Vol. {\bf 9}. 135-162 (1980).


\bibitem{flum} FLUM, J. \& ZIEGLER, M. {\sl Topological model theory.\/} Lecture Notes in Mathematics Vol.
{\bf 769} Springer-Verlag  (1980).

\bibitem{hodges} HODGES, w. {\sl Model theory.\/} Encyclopedia of Mathematics and its Applications. Vol.
{\bf 42} Cambridge Univ. Press (1993).

\bibitem{macintyre} MACINTYRE, A. {\sl Geometrical and Set-Theoretic
    Aspects and Prospects.\/} The Bulletin of Symbolic Logic, Vol. 9,
  No. 2 (Jun., 2003), pp. 197-212

\bibitem{macpherson} MACPHERSON, R.{\sl Intersection Homology and Perverse
Sheaves.\/} {Colloquium Lectures, Annual Meeting of the Amer. Math.
Soc.} San Francisco (1991).



\bibitem{marker} MARKER, D. {\sl Model theory, an introduction}. Graduate Texts in Mathematics. Vol. {\bf 217} Springer-Verlag. (2002). 

\bibitem{mclane} MAC LANE, S. {\sl Categories for the working mathematician}. Graduate Texts in Math. Springer-Verlag.
Vol. {\bf 5}. Springer-Verlag. New York-Heidelberg-Berlin (1971).

\bibitem{mijares} MIJARES, J. \& PADILLA, G. {\sl A Ramsey space of infinite polyhedra and the infinite 
random polyhedron.\/} Arxiv Math.LO {\bf 1209.6421}.

\bibitem{illinois}  SARALEGI, M. {\sl Homological Properties of Stratified Spaces.\/}
Illinois J. Math. {\bf 38}, 47-70 (1994).



\end{thebibliography}
\end{document}